\documentclass[11pt,reqno]{amsart}
\usepackage{amsfonts}
\usepackage{amssymb,latexsym}
\usepackage{cite}
\usepackage[height=190mm,width=138mm]{geometry}

\theoremstyle{plain}
\newtheorem{theorem}{Theorem}

\theoremstyle{definition}

\theoremstyle{remark}

 \numberwithin{equation}{section}
\input{tcilatex}

\begin{document}

\begin{abstract}
In this paper, we use a modular form approach to evaluate the convolution
sums $\sum_{l+42m=n}\sigma (l)\sigma (m)$, $\sum_{2l+21m=n}\sigma (l)\sigma
(m),$ $\sum_{3l+14m=n}\sigma (l)\sigma (m)$ and $\sum_{6l+7m=n}\sigma
(l)\sigma (m) $ for all positive integers $n\mathbf{,}$ and then use their
evaluations to determine the number of representation of \ a positive
integer $n$ by the quadratic form$%
x_{1}^{2}+x_{1}x_{2}+x_{2}^{2}+x_{3}^{2}+x_{3}x_{4}+x_{4}^{2}+14(x_{5}^{2}+x_{5}x_{6}+x_{6}^{2}+x_{7}^{2}+x_{7}x_{8}+x_{8}^{2}).
$\newline

\end{abstract}

\title{Evaluation of the convolution sums $W_{1,42}(n),$ $W_{2,21}(n),$ $%
W_{3,14}(n)$ and $W_{6,7}(n)$}
\author{B\"{u}lent K\"{o}kl\"{u}ce}
\address{The Institute for Computational and Experimental Research in
Mathematics\\
121 S Main St, Providence, 02903\\
Brown University\\
Rhode Island\\
USA} \email{bkokluce@hotmail.com}  \maketitle

\section{Introduction}

Let $\mathbf{%
\mathbb{N}
,}$ $%
\mathbb{Z}
$, $%
\mathbb{Q}
$ and $%
\mathbb{C}
$ denote the set of positive integers, integers, rational numbers and
complex numbers respectively and let $\mathbf{%
\mathbb{N}
}_{\mathbf{0}}=\mathbf{%
\mathbb{N}
\cup \{0\}.}$ For $k\in \mathbf{%
\mathbb{N}
}$ and $n\in
\mathbb{Q}
$ we set

\begin{equation}
\sigma _{k}(n)=\left\{
\begin{array}{cc}
\mathop{\displaystyle \sum }\limits_{_{\substack{ d\in \mathbf{%
\mathbb{N}
}  \\ d\mid n}}}d^{k}, & \text{if }n\in \mathbf{%
\mathbb{N}
,} \\
0, & \text{if }n\in
\mathbb{Q}
,n\notin \mathbf{%
\mathbb{N}
.}%
\end{array}%
\right.  \label{1.1}
\end{equation}
We write\ $\sigma (n)$ for $\sigma _{1}(n)$. Suppose that $r,s\in \mathbf{%
\mathbb{N}
}$ with $r\leq s$ . We define the convolution sum $W_{r,s}(n)$ as follows:

\begin{equation*}
W_{r,s}(n):=\mathop{\displaystyle \sum }\limits_{\substack{ (l,m)\in \mathbf{%
\mathbb{N}
}_{0}^{2}  \\ rl+sm=n}}\sigma (l)\sigma (m).
\end{equation*}

The evaluation of convolution sums $W_{r,s}(n)$ for some levels $rs$
have been done. See Table \ref{Tab:Tb1} for a list of known
convolution sums. In the present
paper, we completed the evaluation of the convolution sums $W_{l,s}(n)$ for $%
(r,s)=(1,42),(2,21),(3,14)$ and $(6,7)$ by using a modular form approach. We
also evaluated the results for $(1,14),(2,7),(1,7)$ which were firstly given
in \cite{Alaca10} and see that they are consistent.

\begin{table}[th]
\caption{A list of previously known convolution sums.}
\label{Tab:Tb1}\centering

\begin{tabular}{|l|l|l|}
\hline
\ \ Level $rs$ & \ \ \ \ \ \ \ \ \ \ \ Authors & \ \ \ \ \ \ \ References \\
\hline\hline
$1$ & $%
\begin{array}{c}
\text{M.Besge \& J.W.L. Glaisher \&} \\
\text{S. Ramanujan}%
\end{array}%
$ & \cite{5Besge,Glaisher,Ramanujan} \\ \hline
$2,3,4$ & $%
\begin{array}{c}
\text{J. G. Huard \& Z. M. Ou \&} \\
\text{ B. K. Spearman \& K. S. Williams}%
\end{array}%
$ & \cite{Huard} \\ \hline
$5,7$ & $%
\begin{array}{c}
\text{M. Lemire \& K. S. Williams \&} \\
\text{S. Cooper \& P. C. Toh}%
\end{array}%
$ & \cite{Lemire,Cooper2} \\ \hline
$6$ & \c{S}. Alaca \& K. S. Williams & \cite{Alaca6} \\ \hline
$8,9$ & K. S. Willams & \cite{Williams3,Williams4} \\ \hline
$10,11,13,14$ & E. Royer & \cite{Royer} \\ \hline
$12,16,18,24$ & A. Alaca \& \c{S}. Alaca \& K. S. Williams & \cite%
{Alaca1,Alaca3,Alaca2,Alaca8} \\ \hline
$15$ & B. Ramakrishnan \& B. Sabu & \cite{Ramakrishnan2} \\ \hline
$10,20$ & S. Cooper \& D. Ye & \cite{Cooper} \\ \hline
$23$ & H. H. Chan \& S. Cooper & \cite{Chan} \\ \hline
$25$ & $%
\begin{array}{c}
\text{E. X. W. Xia \& X. L. Tian \&} \\
\text{O. X. M. Yao}%
\end{array}%
$ & \cite{Xia} \\ \hline
$27,32,48,54$ & \c{S}. Alaca \& Y. Kesicio\u{g}lu & \cite{Alaca7}, \cite%
{Alacay} \\ \hline $36$ & D. Ye & \cite{Ye} \\ \hline $14,26,28,30$
& A. Alaca \& \c{S}. Alaca \& E. Ntienjem & \cite{Alaca10},
\cite{Ebe1} \\ \hline $22,44,52$ & E. Ntienjem & \cite{Ebe3} \\
\hline $27,40,55$ & B. Kendirli & \cite{Ken} \\ \hline $33,40,56$ &
E. Ntienjem & \cite{Ebe4} \\ \hline $48,64$ & E. Ntienjem &
\cite{Ebe2} \\ \hline $17,34,68$ & B. K\"{o}kl\"{u}ce & \cite{Kokluce3} \\
\hline
\end{tabular}%
\end{table}

For $l\in \mathbf{%
\mathbb{N}
}$ and $n\in \mathbf{%
\mathbb{N}
}_{\mathbf{0}}$ we let $N_{l}(n)$ denote the representation number of $n$ by
the form $%
x_{1}^{2}+x_{1}x_{2}+x_{2}^{2}+x_{3}^{2}+x_{3}x_{4}+x_{4}^{2}+l(x_{5}^{2}+x_{5}x_{6}+x_{6}^{2}+x_{7}^{2}+x_{7}x_{8}+x_{8}^{2}),
$ that is
\begin{equation}
N_{l}(n):=\text{card}\left\{
\begin{array}{c}
\left( x_{1},...,x_{8}\right) \in
\mathbb{Z}
^{8}:n=x_{1}^{2}+x_{1}x_{2}+x_{2}^{2}+x_{3}^{2}+x_{3}x_{4}+x_{4}^{2} \\
+l(x_{5}^{2}+x_{5}x_{6}+x_{6}^{2}+x_{7}^{2}+x_{7}x_{8}+x_{8}^{2})%
\end{array}%
\right\} \text{.}  \label{1.2}
\end{equation}
Explicit formulae for $N_{l}(n)$ are obtained before for $l\leq 12$ and for $%
l=16$ and $18$. See for example, \cite%
{Lomadze,Alaca6,Williams3,Alaca1,Alaca2,Alaca3,Alaca10,Alaca7,Alaca8,Ramakrishnan2,Ebe1,Alacay}%
. The author of this article also found many such representation number
formulae for the direct sum of quadratic forms of this type, see for
example, \cite{Kokluce1,Kokluce2}. In this article, we use the convolutions
sums $W_{3,14}(n)\ $and $W_{1,42}(n)$ obtained here with the convolution sum
$W_{1,14}(n)$ to find explicit formula for $N_{14}(n)$.

The rest of this paper is organized as follows. In Section 2, we give some
preliminary results related to the eta products and modular forms$.$ In
Section 3, we give our main theorem about the convolution sums of divisor
functions and its proof. In Section 4, we give a formula for $N_{14}(n)$ and
then prove it.

\section{Preliminary Results}

For $N\in \mathbf{%
\mathbb{N}
}$ and $k\in
\mathbb{Z}
$ we write $M_{k}(\Gamma _{0}(N))$ to denote the space of modular forms of
weight $k$ (with trivial multiplier system) for the modular subgroup $\Gamma
_{0}(N)$ defined by

\begin{equation}
\Gamma _{0}(N)=\left\{ \left(
\begin{array}{cc}
a & b \\
c & d%
\end{array}%
\right) \in SL_{2}(%
\mathbb{Z}
):c\equiv 0(\func{mod}N)\right\} .  \label{1.3}
\end{equation}
It is known (see, for example \cite[p.83]{Stein}) that
\begin{equation}
M_{k}(\Gamma _{0}(N))=E_{k}(\Gamma _{0}(N))\oplus S_{k}(\Gamma _{0}(N))
\label{e0}
\end{equation}
where $E_{k}(\Gamma _{0}(N))$ and $S_{k}(\Gamma _{0}(N))$ are the
corresponding subspaces of Eisenstein forms and cusp forms of weight $k$ for
the modular subgroup $\Gamma _{0}(N).$

The Dedekind eta function $\eta (z)$ is the holomorphic function defined on
the upper half plane $\mathbb{H}=\left\{ z\in
\mathbb{C}
:\func{Im}(z)>0\right\} $ by the product formula

\begin{equation}
\eta (z)=e^{\pi iz/12}\mathop{\displaystyle \prod }\limits_{n=1}^{\infty
}(1-e^{2\pi inz}).  \label{kk}
\end{equation}
Through the remainder of the paper we take $q=q(z):=e^{2\pi iz}$ with $z\in
\mathbb{H}$ and so by (\ref{kk}) we have

\begin{equation}
\eta (z)=q^{1/24}\mathop{\displaystyle \prod }\limits_{n=1}^{\infty
}(1-q^{n}).  \label{e1}
\end{equation}
An eta quotient is defined to be a finite product of the form

\begin{equation}
f(z)=\mathop{\displaystyle \prod }\limits_{\delta }\eta ^{r_{\delta
}}(\delta z),  \label{e2}
\end{equation}
where $\delta $ runs through a finite set of positive integers and the
exponents $r_{\delta }$ are nonzero integers. By taking $N$ to be the least
common multiple of $\delta ^{\prime }$s we can write the eta quotient (\ref%
{e2}) as

\begin{equation}
f(z)=\mathop{\displaystyle \prod }\limits_{\delta \mid N}\eta ^{r_{\delta
}}(\delta z)  \label{e3}
\end{equation}
where some of the exponents $r_{\delta }$ may be zero. When all of the
exponents $r_{\delta }$ are nonnegative, $f(z)$ is said to be an eta
product. Now we define the following $20$ eta quotients

\begin{equation}
C_{1}(q)=\frac{\eta ^{5}(z)\eta ^{5}(7z)}{\eta (2z)\eta (14z)},  \label{t1}
\end{equation}

\begin{equation}
C_{2}(q)=\eta ^{2}(z)\eta ^{2}(2z)\eta ^{2}(7z)\eta ^{2}(14z),  \label{t2}
\end{equation}

\begin{equation}
C_{3}(q)=\frac{\eta ^{6}(z)\eta ^{6}(14z)}{\eta ^{2}(2z)\eta ^{2}(7z)},
\label{t3}
\end{equation}

\begin{equation}
C_{4}(q)=\frac{\eta ^{6}(2z)\eta ^{6}(7z)}{\eta ^{2}(z)\eta ^{2}(14z)},
\label{t4}
\end{equation}

\begin{equation}
C_{5}(q)=\eta ^{2}(z)\eta ^{2}(2z)\eta ^{2}(3z)\eta ^{2}(6z),
\end{equation}

\begin{equation}
C_{6}(q)=\eta (2z)\eta ^{3}(3z)\eta (7z)\eta ^{3}(42z),
\end{equation}

\begin{equation}
C_{7}(q)=\eta (z)\eta ^{3}(6z)\eta (14z)\eta ^{3}(21z),
\end{equation}

\begin{equation}
C_{8}(q)=\eta ^{2}(7z)\eta ^{2}(14z)\eta ^{2}(21z)\eta ^{2}(42z),
\end{equation}

\begin{equation}
C_{9}(q)=\eta (3z)\eta (6z)\eta ^{3}(7z)\eta ^{3}(14z),
\end{equation}

\begin{equation}
C_{10}(q)=\eta ^{2}(2z)\eta ^{2}(6z)\eta ^{2}(7z)\eta ^{2}(21z),
\end{equation}

\begin{equation}
C_{11}(q)=\eta (z)\eta (2z)\eta ^{3}(21z)\eta ^{3}(42z),
\end{equation}

\begin{equation}
C_{12}(q)=\eta ^{3}(z)\eta (6z)\eta ^{3}(14z)\eta (21z),
\end{equation}

\begin{equation}
C_{13}(q)=\eta ^{3}(z)\eta ^{3}(2z)\eta (21z)\eta (42z),
\end{equation}

\begin{equation}
C_{14}(q)=\frac{\eta ^{2}(2z)\eta ^{2}(14z)\eta ^{6}(21z)}{\eta ^{2}(7z)},
\end{equation}

\begin{equation}
C_{15}(q)=\frac{\eta ^{6}(2z)\eta ^{2}(3z)\eta ^{2}(21z)}{\eta ^{2}(6z)},
\end{equation}

\begin{equation}
C_{16}(q)=\frac{\eta ^{5}(2z)\eta ^{5}(21z)}{\eta (z)\eta (42z)},
\end{equation}

\begin{equation}
C_{17}(q)=\frac{\eta ^{6}(6z)\eta ^{6}(21z)}{\eta ^{2}(3z)\eta ^{2}(42z)},
\end{equation}

\begin{equation}
C_{18}(q)=\frac{\eta ^{3}(6z)\eta (7z)\eta (14z)\eta ^{4}(21z),}{\eta (3z)}
\end{equation}

\begin{equation}
C_{19}(q)=\eta ^{2}(3z)\eta ^{2}(6z)\eta ^{2}(21z)\eta ^{2}(42z),
\end{equation}

\begin{equation}
C_{20}(q)=\frac{\eta ^{4}(3z)\eta ^{4}(14z)\eta ^{2}(21z)\eta ^{2}(42z)}{%
\eta ^{2}(6z)\eta ^{2}(7z)},  \label{t20}
\end{equation}
and the integers $c_{k}(n)$ $(n\in
\mathbb{N}
)$ for $1\leq k\leq 20$ by

\begin{equation}
C_{k}(q)=\mathop{\displaystyle \sum }\limits_{n=1}^{\infty }c_{k}(n)q^{n}.
\label{ck}
\end{equation}
The Eisenstein series $L(q)$ and $M(q)$ are defined by

\begin{equation}
L(q)=E_{2}(q)=1-24\mathop{\displaystyle \sum }\limits_{n=1}^{\infty }\sigma
\left( n\right) q^{n},  \label{2.1}
\end{equation}
and
\begin{equation}
M(q)=E_{4}(q)=1+240\mathop{\displaystyle \sum }\limits_{n=1}^{\infty }\sigma
_{3}\left( n\right) q^{n}.  \label{2.2}
\end{equation}
We use the following theorem to determine if a given eta product is in $%
M_{k}(\Gamma _{0}(N)),$ see \cite[Theorem 5.7, p.99]{Kilford}. Note that the
eta quotients given in (\ref{t1})-(\ref{t4}) are the same\ first four eta
quotients which are used in \cite{Alaca10}.

\begin{theorem}
\label{thm:1} Let $N$ be a positive integer and let $f(z)=%
\mathop{\displaystyle \prod }\nolimits_{1\leq \delta \mid N}\eta ^{r_{\delta
}}(\delta z)$ be an eta quotient which satisfies the following conditions:

\begin{enumerate}
\item[(i)] $\mathop{\displaystyle \sum }\nolimits_{1\leq \delta \mid
N}\delta .r_{\delta }\equiv 0(\func{mod}24),$

\item[(ii)] $\mathop{\displaystyle \sum }\nolimits_{1\leq \delta \mid N}%
\dfrac{N}{\delta }.r_{\delta }\equiv 0(\func{mod}24),$

\item[(iii)] $\mathop{\displaystyle \prod }\nolimits_{1\leq \delta \mid
N}\delta ^{r_{\delta }}$ is the square of a rational number,

\item[(iv)] for each $d\mid N,\mathop{\displaystyle \sum }\nolimits_{1\leq
\delta \mid N}\dfrac{\gcd (d,\delta )^{2}.r_{\delta }}{\delta }\geq 0,$

\item[(v)] the weight $k=\dfrac{1}{2}\mathop{\displaystyle \sum }%
\nolimits_{1\leq \delta \mid N}r_{\delta }$ is an even integer.

Then $f(z)$ is in $M_{k}(\Gamma _{0}(N)).$ In addition to the above
conditions if all the inequalities in (iv) hold strictly then $f(z)$ is in $%
S_{k}(\Gamma _{0}(N)).$
\end{enumerate}
\end{theorem}

Note that the cusp forms $C_{k}(q)$ $(1\leq k\leq 20)$ defined in (\ref{t1}%
)-(\ref{t20}) are constructed in a way that they satisfy the conditions of
Theorem \ref{thm:1} from (i) to (v).

\begin{theorem}
\label{thm:2} (a) $\left\{ C_{k}(q)\text{ }(1\leq k\leq 20)\right\} $ is a
basis for $S_{4}(\Gamma _{0}(42)).$

(b) $E_{4}(q^{t})$ $(t=1,2,3,6,7,14,21,42)$ constitute a basis for $%
E_{4}(\Gamma _{0}(42)).$

(c) $\left\{ C_{k}(q)(1\leq k\leq 20)\right\} $ together with $E_{4}(q^{t})$
$(t=1,2,3,6,7,14,21,42)$ constitute a basis for $M_{4}(\Gamma _{0}(42)).$
\end{theorem}

\begin{proof}
(a) It follows from (\ref{t1})-(\ref{t20}), (\ref{e1})-(\ref{e3})
and Theorem \ref{thm:1} that $C_{k}(q)(1\leq k\leq 20)$ are in
$S_{4}(\Gamma _{0}(42)).$ By \cite[Theorem 3.8, p.50]{Kilford}, the
dimension of $S_{4}(\Gamma _{0}(42))$ is $20$. We use the Maple
software to show that there is no linear relationship among
$C_{k}(q)(1\leq k\leq 20)$. Thus $C_{k}(q)(1\leq k\leq 20)$ form a
basis of $S_{4}(\Gamma _{0}(42)).$

(b) It can be shown by using the dimension formula for the Eisenstein space
(see for example \cite[Theorem 3.8, p.50]{Kilford}) that the dimension of $%
E_{4}(\Gamma _{0}(42))$ is $8$. Thus it follows from \cite[Theorem 5.9, p.88]%
{Stein} that $E_{4}(q^{t})$ $(t=1,2,3,6,7,14,21,42)$ constitute a basis of $%
E_{4}(\Gamma _{0}(42)).$

(c) It follows from parts (a), (b) and (\ref{e0}) that the dimension of $%
M_{4}(\Gamma _{0}(42))$ is $28$, and therefore $E_{4}(q^{t})$ $%
(t=1,2,3,6,7,14,21,42)$ together with $C_{k}(q)$ $(1\leq k\leq 20)$
constitute a basis for the space $M_{4}(\Gamma _{0}(42)).$
\end{proof}

The following theorem is given in \cite[Proposition 2.12, p 23]{Kilford}.

\begin{theorem}
\label{thm:3} Let $N$ be a positive integer. Then

\begin{equation*}
M=\left[ SL_{2}(%
\mathbb{Z}
):\Gamma _{0}(N)\right] =N\mathop{\displaystyle \prod }\limits_{_{p\mid
N}}(1+\frac{1}{p}).
\end{equation*}
\end{theorem}

We use Theorem \ref{thm:3} with Sturm's bound theorem (see\cite[Theorem
3.13, p.53]{Kilford} or \cite{Sturm}) to see if two modular forms in the
same space are equal. We can restate the Sturm's theorem for our case as
follows.

\begin{theorem}
\label{thm:4} Let $\Gamma _{0}\in SL_{2}(%
\mathbb{Z}
)$ be a congruence subgroup of index $M$ and let $f\in M_{4}(\Gamma _{0}(N))$
be a modular form. If%
\begin{equation*}
\upsilon _{\infty }(f)>\frac{M}{3}=S(N)
\end{equation*}
then $f$ is identically zero. Thus, if $f_{1}(z)$ and $f_{2}(z)$ are two
modular forms in $M_{4}(\Gamma _{0}(N))$ with Fourier series expansions $%
f_{1}(z)=\mathop{\displaystyle \sum }\limits_{n=1}^{\infty }a_{n}q^{n}$ and $%
f_{2}(z)=\mathop{\displaystyle \sum }\limits_{n=1}^{\infty }b_{n}q^{n}$ such
that $a_{n}=b_{n}$ for all $n\leq S(N)$ then $f(z)=g(z)$. By Theorem 4, it
is clear that the Sturm bound for $M_{4}(\Gamma _{0}(42))$ is
\end{theorem}

\begin{equation}
S(42)=32.  \label{st}
\end{equation}

\section{Evaluation of Convolution Sums}

\begin{theorem}
\label{thm:5} \bigskip We have,%
\begin{eqnarray*}
(L(q)-42L(q^{42}))^{2} &=&\frac{604}{625}M(q)-\frac{84}{625}M(q^{2})-\frac{%
189}{625}M(q^{3})-\frac{756}{625}M(q^{6}) \\
&&-\frac{1029}{625}M(q^{7})-\frac{4116}{625}M(q^{14})-\frac{9261}{625}%
M(q^{21}) \\
&&+\frac{1065456}{625}M(q^{42})+\frac{6912}{5}C_{1}(q)+\frac{24624}{25}%
C_{2}(q) \\
&&+\frac{1728}{5}C_{4}(q)+\frac{1008}{125}C_{5}(q)-\frac{296352}{5}C_{6}(q)
\\
&&+\frac{5340384}{25}C_{7}(q)-\frac{2067408}{125}C_{8}(q)-\frac{34272}{5}%
C_{9}(q) \\
&&+\frac{2346624}{25}C_{10}(q)-\frac{3284064}{5}C_{11}(q)+\frac{1024128}{25}%
C_{12}(q) \\
&&-\frac{124992}{5}C_{13}(q)+\frac{653184}{5}C_{14}(q)+\frac{48384}{5}%
C_{15}(q) \\
&&-\frac{520128}{5}C_{16}(q)+\frac{217728}{5}C_{18}(q)+\frac{27216}{25}%
C_{19}(q) \\
&&+\frac{36288}{5}C_{20}(q)
\end{eqnarray*}

\begin{eqnarray*}
(2L(q^{2})-21L(q^{21}))^{2} &=&-\frac{21}{625}M(q)+\frac{2416}{625}M(q^{2})-%
\frac{189}{625}M(q^{3})-\frac{756}{625}M(q^{6}) \\
&&-\frac{1029}{625}M(q^{7})-\frac{4116}{625}M(q^{14})+\frac{266364}{625}%
M(q^{21}) \\
&&-\frac{37044}{625}M(q^{42})+\frac{3456}{5}C_{1}(q)+\frac{24624}{25}C_{2}(q)
\\
&&+1728C_{3}(q)-\frac{3456}{5}C_{4}(q)+\frac{1008}{125}C_{5}(q) \\
&&-78624C_{6}(q)+\frac{4251744}{25}C_{7}(q)+\frac{4282992}{125}C_{8}(q) \\
&&-\frac{82656}{5}C_{9}(q)+\frac{2346624}{25}C_{10}(q)-\frac{2292192}{5}%
C_{11}(q) \\
&&+\frac{983808}{25}C_{12}(q)-\frac{133056}{5}C_{13}(q)+\frac{870912}{5}%
C_{14}(q) \\
&&+\frac{24192}{5}C_{15}(q)-\frac{471744}{5}C_{16}(q)-\frac{72576}{5}%
C_{18}(q) \\
&&+\frac{27216}{25}C_{19}(q)+\frac{108864}{5}C_{20}(q)
\end{eqnarray*}

\begin{eqnarray*}
(3L(q^{3})-14L(q^{14}))^{2} &=&-\frac{21}{625}M(q)-\frac{84}{625}M(q^{2})+%
\frac{5436}{625}M(q^{3})-\frac{756}{625}M(q^{6}) \\
&&-\frac{1029}{625}M(q^{7})+\frac{118384}{625}M(q^{14})-\frac{9261}{625}%
M(q^{21}) \\
&&-\frac{37044}{625}M(q^{42})-\frac{6912}{5}C_{1}(q)+\frac{3024}{25}C_{2}(q)
\\
&&+\frac{6912}{5}C_{4}(q)+\frac{1008}{125}C_{5}(q)+\frac{223776}{5}C_{6}(q)
\\
&&-\frac{5062176}{25}C_{7}(q)+\frac{2166192}{125}C_{8}(q)+\frac{34272}{5}%
C_{9}(q) \\
&&-\frac{2491776}{25}C_{10}(q)+\frac{3284064}{5}C_{11}(q) \\
&&-\frac{1052352}{25}C_{12}(q)+\frac{124992}{5}C_{13}(q)-\frac{653184}{5}%
C_{14}(q) \\
&&-\frac{48384}{5}C_{15}(q)+\frac{520128}{5}C_{16}(q)+15552C_{17}(q) \\
&&-\frac{217728}{5}C_{18}(q)+\frac{221616}{25}C_{19}(q)-\frac{36288}{5}%
C_{20}(q)
\end{eqnarray*}

\begin{eqnarray*}
(6L(q^{6})-7L(q^{7}))^{2} &=&-\frac{21}{625}M(q)-\frac{84}{625}M(q^{2})-%
\frac{189}{625}M(q^{3})+\frac{21744}{625}M(q^{6}) \\
&&+\frac{29596}{625}M(q^{7})-\frac{4116}{625}M(q^{14})-\frac{9261}{625}%
M(q^{21}) \\
&&-\frac{37044}{625}M(q^{42})-\frac{3456}{5}C_{1}(q)+\frac{175824}{25}%
C_{2}(q) \\
&&-1728C_{3}(q)-\frac{5184}{5}C_{4}(q)+\frac{217008}{125}C_{5}(q) \\
&&+\frac{320544}{5}C_{6}(q)-\frac{5960736}{25}C_{7}(q)+\frac{6399792}{125}%
C_{8}(q) \\
&&+\frac{82656}{5}C_{9}(q)-\frac{763776}{25}C_{10}(q)+\frac{2292192}{5}%
C_{11}(q) \\
&&-\frac{1674432}{25}C_{12}(q)+\frac{133056}{5}C_{13}(q)-\frac{870912}{5}%
C_{14}(q) \\
&&-\frac{24192}{5}C_{15}(q)+\frac{471744}{5}C_{16}(q)-15552C_{17}(q) \\
&&+\frac{72576}{5}C_{18}(q)+\frac{1776816}{25}C_{19}(q)-\frac{108864}{5}%
C_{20}(q)
\end{eqnarray*}
\end{theorem}

\begin{proof}
We only prove the case for $(L(q)-42L(q^{42}))^{2}$ as the rest can be
proven in a similar way. By \cite[Theorem 5.8]{Stein} we have $%
L(q)-42L(q^{42})\in M_{2}(\Gamma _{0}(42))$, and so

\begin{equation*}
(L(q)-42L(q^{42}))^{2}\in M_{4}(\Gamma _{0}(42)).
\end{equation*}
Thus, appealing to Theorem \ref{thm:1}(c), we can express $%
(L(q)-42L(q^{42}))^{2}$ as a linear combination of $E_{4}(q^{t})$ $%
(t=1,2,3,6,7,14,21,42)$ and $C_{k}(q)(1\leq k\leq 20)$. So, there exist
coefficients $x_{t}(t\in
\mathbb{N}
$, $t\mid 42)$ and $y_{k}(k=1,...,20)$ such that

\begin{equation}
(L(q)-42L(q^{42}))^{2}=\mathop{\displaystyle \sum }\limits_{t\in
\mathbb{N}
,t\mid 42}x_{t}E_{4}(q^{t})+\mathop{\displaystyle \sum }%
\limits_{k=1}^{20}y_{k}C_{k}(q).  \label{lin}
\end{equation}
Appealing to (\ref{st}) and equating the coefficients of $q^{n}$ for $1\leq
n\leq 32$ on both sides of (\ref{lin}) we obtain required result. The
following theorem can be given as the result of Theorem \ref{thm:5},
equations (\ref{ck}) and (\ref{2.2}).
\end{proof}

\begin{theorem}
\label{thm:6}
\begin{eqnarray*}
(L(q)-42L(q^{42}))^{2} &=&1681+\sum_{n=1}^{\infty }(\frac{28992}{125}\sigma
_{3}(n)-\frac{4032}{125}\sigma _{3}(\frac{n}{2})-\frac{9072}{125}\sigma _{3}(%
\frac{n}{3}) \\
&&-\frac{36288}{125}\sigma _{3}(\frac{n}{6})-\frac{49392}{125}\sigma _{3}(%
\frac{n}{7})-\frac{197568}{125}\sigma _{3}(\frac{n}{14}) \\
&&-\frac{444528}{125}\sigma _{3}(\frac{n}{21})+\frac{51141888}{125}\sigma
_{3}(\frac{n}{42})+\frac{6912}{5}c_{1}(n) \\
&&+\frac{24624}{25}c_{2}(n)+\frac{1728}{5}c_{4}(n)+\frac{1008}{125}c_{5}(n)-%
\frac{296352}{5}c_{6}(n) \\
&&+\frac{5340384}{25}c_{7}(n)-\frac{2067408}{125}c_{8}(n)-\frac{34272}{5}%
c_{9}(n) \\
&&+\frac{2346624}{25}c_{10}(n)-\frac{3284064}{5}c_{11}(n)+\frac{1024128}{25}%
c_{12}(n) \\
&&-\frac{124992}{5}c_{13}(n)+\frac{653184}{5}c_{14}(n)+\frac{48384}{5}%
c_{15}(n) \\
&&-\frac{520128}{5}c_{16}(n)+\frac{217728}{5}c_{18}(n)+\frac{27216}{25}%
c_{19}(n) \\
&&+\frac{36288}{5}c_{20}(n))q^{n}
\end{eqnarray*}

\begin{eqnarray*}
(2L(q^{2})-21L(q^{21}))^{2} &=&361+\sum_{n=1}^{\infty }(-\frac{1008}{125}%
\sigma _{3}(n)+\frac{115968}{125}\sigma _{3}(\frac{n}{2})-\frac{9072}{125}%
\sigma _{3}(\frac{n}{3}) \\
&&-\frac{36288}{125}\sigma _{3}(\frac{n}{6})-\frac{49392}{125}\sigma _{3}(%
\frac{n}{7})-\frac{197568}{125}\sigma _{3}(\frac{n}{14}) \\
&&+\frac{12785472}{125}\sigma _{3}(\frac{n}{21})-\frac{1778112}{125}\sigma
_{3}(\frac{n}{42})+\frac{3456}{5}c_{1}(n) \\
&&+\frac{24624}{25}c_{2}(n)+1728c_{3}(n)-\frac{3456}{5}c_{4}(n)+\frac{1008}{%
125}c_{5}(n) \\
&&-78624c_{6}(n)+\frac{4251744}{25}c_{7}(n)+\frac{4282992}{125}c_{8}(n) \\
&&-\frac{82656}{5}c_{9}(n)+\frac{2346624}{25}c_{10}(n)-\frac{2292192}{5}%
c_{11}(n) \\
&&+\frac{983808}{25}c_{12}(n)-\frac{133056}{5}c_{13}(n)+\frac{870912}{5}%
c_{14}(n) \\
&&+\frac{24192}{5}c_{15}(n)-\frac{471744}{5}c_{16}(n)-\frac{72576}{5}%
c_{18}(n) \\
&&+\frac{27216}{25}c_{19}(n)+\frac{108864}{5}c_{20}(n))q^{n}
\end{eqnarray*}

\begin{eqnarray*}
(3L(q^{3})-14L(q^{14}))^{2} &=&121+\sum_{n=1}^{\infty }(-\frac{1008}{125}%
\sigma _{3}(n)-\frac{4032}{125}\sigma _{3}(\frac{n}{2})+\frac{260928}{125}%
\sigma _{3}(\frac{n}{3}) \\
&&-\frac{36288}{125}\sigma _{3}(\frac{n}{6})-\frac{49392}{125}\sigma _{3}(%
\frac{n}{7})+\frac{5682432}{125}\sigma _{3}(\frac{n}{14}) \\
&&-\frac{444528}{125}\sigma _{3}(\frac{n}{21})-\frac{1778112}{125}\sigma
_{3}(\frac{n}{42})-\frac{6912}{5}c_{1}(n) \\
&&+\frac{3024}{25}c_{2}(n)+\frac{6912}{5}c_{4}(n)+\frac{1008}{125}c_{5}(n) \\
&&+\frac{223776}{5}c_{6}(n)-\frac{5062176}{25}c_{7}(n)+\frac{2166192}{125}%
c_{8}(n) \\
&&+\frac{34272}{5}c_{9}(n)-\frac{2491776}{25}c_{10}(n)+\frac{3284064}{5}%
c_{11}(n) \\
&&-\frac{1052352}{25}c_{12}(n)+\frac{124992}{5}c_{13}(n)-\frac{653184}{5}%
c_{14}(n) \\
&&-\frac{48384}{5}c_{15}(n)+\frac{520128}{5}c_{16}(n)+15552c_{17}(n) \\
&&-\frac{217728}{5}c_{18}(n)+\frac{221616}{25}c_{19}(n)-\frac{36288}{5}%
c_{20}(n))q^{n}
\end{eqnarray*}
\end{theorem}

\begin{eqnarray*}
(6L(q^{6})-7L(q^{7}))^{2} &=&1+\sum_{n=1}^{\infty }(-\frac{1008}{125}\sigma
_{3}(n)-\frac{4032}{125}\sigma _{3}(\frac{n}{2})-\frac{9072}{125}\sigma _{3}(%
\frac{n}{3}) \\
&&+\frac{1043712}{125}\sigma _{3}(\frac{n}{6})+\frac{1420608}{125}\sigma
_{3}(\frac{n}{7})-\frac{197568}{125}\sigma _{3}(\frac{n}{14}) \\
&&-\frac{444528}{125}\sigma _{3}(\frac{n}{21}))-\frac{1778112}{125}\sigma
_{3}(\frac{n}{42})-\frac{3456}{5}c_{1}(n) \\
&&+\frac{175824}{25}c_{2}(n)-1728c_{3}-\frac{5184}{5}c_{4}+\frac{217008}{125}%
c_{5} \\
&&+\frac{320544}{5}c_{6}(n)-\frac{5960736}{25}c_{7}(n)+\frac{6399792}{125}%
c_{8}(n) \\
&&+\frac{82656}{5}c_{9}(n)-\frac{763776}{25}c_{10}(n)+\frac{2292192}{5}%
c_{11}(n) \\
&&-\frac{1674432}{25}c_{12}(n)+\frac{133056}{5}c_{13}(n)-\frac{870912}{5}%
c_{14}(n) \\
&&-\frac{24192}{5}c_{15}(n)+\frac{471744}{5}c_{16}(n)-15552c_{17}(n) \\
&&+\frac{72576}{5}c_{18}(n)+\frac{1776816}{25}c_{19}(n)-\frac{108864}{5}%
c_{20}(n))q^{n}
\end{eqnarray*}

\begin{theorem}
\label{thm:7} Let $n\in
\mathbb{N}
.$ Then
\begin{eqnarray}
W_{1,42}(n) &=&\frac{1}{6000}\sigma _{3}\left( n\right) +\frac{1}{1500}%
\sigma _{3}\left( \frac{n}{2}\right) +\frac{3}{2000}\sigma _{3}\left( \frac{n%
}{3}\right) +\frac{3}{500}\sigma _{3}\left( \frac{n}{6}\right)  \label{w1} \\
&&+\frac{49}{6000}\sigma _{3}\left( \frac{n}{7}\right) +\frac{49}{1500}%
\sigma _{3}\left( \frac{n}{14}\right) +\frac{147}{2000}\sigma _{3}\left(
\frac{n}{21}\right)  \notag \\
&&+\frac{147}{500}\sigma _{3}\left( \frac{n}{42}\right) +(\frac{1}{24}-\frac{%
n}{168})\sigma \left( n\right) +(\frac{1}{24}-\frac{n}{4})\sigma \left(
\frac{n}{42}\right)  \notag \\
&&-\frac{1}{35}c_{1}(n)-\frac{57}{2800}c_{2}(n)-\frac{1}{140}c_{4}(n)-\frac{1%
}{6000}c_{5}(n)  \notag \\
&&+\frac{49}{40}c_{6}(n)-\frac{883}{200}c_{7}(n)+\frac{2051}{6000}c_{8}(n)+%
\frac{17}{120}c_{9}(n)-\frac{97}{50}c_{10}(n)  \notag \\
&&+\frac{543}{40}c_{11}(n)-\frac{127}{150}c_{12}(n)+\frac{31}{60}c_{13}(n)-%
\frac{27}{10}c_{14}(n)-\frac{1}{5}c_{15}(n)  \notag \\
&&+\frac{43}{20}c_{16}(n)-\frac{9}{10}c_{18}(n)-\frac{9}{400}c_{19}(n)-\frac{%
3}{20}c_{20}(n)  \notag
\end{eqnarray}

\begin{eqnarray}
W_{2,21}(n) &=&\frac{1}{6000}\sigma _{3}\left( n\right) +\frac{1}{1500}%
\sigma _{3}\left( \frac{n}{2}\right) +\frac{3}{2000}\sigma _{3}\left( \frac{n%
}{3}\right) +\frac{3}{500}\sigma _{3}\left( \frac{n}{6}\right)  \label{w2} \\
&&+\frac{49}{6000}\sigma _{3}\left( \frac{n}{7}\right) +\frac{49}{1500}%
\sigma _{3}\left( \frac{n}{14}\right) +\frac{147}{2000}\sigma _{3}\left(
\frac{n}{21}\right)  \notag \\
&&+\frac{147}{500}\sigma _{3}\left( \frac{n}{42}\right) +(\frac{1}{24}-\frac{%
n}{84})\sigma \left( \frac{n}{2}\right) +(\frac{1}{24}-\frac{n}{8})\sigma
\left( \frac{n}{21}\right)  \notag \\
&&-\frac{1}{70}c_{1}(n)-\frac{57}{2800}c_{2}(n)-\frac{1}{28}c_{3}(n)+\frac{1%
}{70}c_{4}(n)-\frac{1}{6000}c_{5}(n)  \notag \\
&&+\frac{13}{8}c_{6}(n)-\frac{703}{200}c_{7}(n)-\frac{4249}{6000}c_{8}(n)+%
\frac{41}{120}c_{9}(n)-\frac{97}{50}c_{10}(n)  \notag \\
&&+\frac{379}{40}c_{11}(n)-\frac{61}{75}c_{12}(n)+\frac{11}{20}c_{13}(n)-%
\frac{18}{5}c_{14}(n)-\frac{1}{10}c_{15}(n)  \notag \\
&&+\frac{39}{20}c_{16}(n)+\frac{3}{10}c_{18}(n)-\frac{9}{400}c_{19}(n)-\frac{%
9}{20}c_{20}(n)  \notag
\end{eqnarray}

\begin{eqnarray}
W_{3,14}(n) &=&\bigskip \frac{1}{6000}\sigma _{3}\left( n\right) +\frac{1}{%
1500}\sigma _{3}\left( \frac{n}{2}\right) +\frac{3}{2000}\sigma _{3}\left(
\frac{n}{3}\right) +\frac{3}{500}\sigma _{3}\left( \frac{n}{6}\right)
\label{w3} \\
&&+\frac{49}{6000}\sigma _{3}\left( \frac{n}{7}\right) +\frac{49}{1500}%
\sigma _{3}\left( \frac{n}{14}\right) +\frac{147}{2000}\sigma _{3}\left(
\frac{n}{21}\right)  \notag \\
&&+\frac{147}{500}\sigma _{3}\left( \frac{n}{42}\right) +(\frac{1}{24}-\frac{%
n}{56})\sigma \left( \frac{n}{3}\right) +(\frac{1}{24}-\frac{n}{12})\sigma
\left( \frac{n}{14}\right)  \notag \\
&&+\frac{1}{35}c_{1}(n)-\frac{1}{400}c_{2}(n)-\frac{1}{35}c_{4}(n)-\frac{1}{%
6000}c_{5}(n)-\frac{37}{40}c_{6}(n)  \notag \\
&&+\frac{837}{200}c_{7}(n)-\frac{2149}{6000}c_{8}(n)-\frac{17}{120}c_{9}(n)+%
\frac{103}{50}c_{10}(n)  \notag \\
&&-\frac{543}{40}c_{11}(n)+\frac{87}{100}c_{12}(n)-\frac{31}{60}c_{13}(n)+%
\frac{27}{10}c_{14}(n)+\frac{1}{5}c_{15}(n)  \notag \\
&&-\frac{43}{20}c_{16}(n)-\frac{9}{28}c_{17}(n)+\frac{9}{10}c_{18}(n)-\frac{%
513}{2800}c_{19}(n)+\frac{3}{20}c_{20}(n)  \notag
\end{eqnarray}

\begin{eqnarray}
W_{6,7}(n) &=&\frac{1}{6000}\sigma _{3}\left( n\right) +\frac{1}{1500}\sigma
_{3}\left( \frac{n}{2}\right) +\frac{3}{2000}\sigma _{3}\left( \frac{n}{3}%
\right) +\frac{3}{500}\sigma _{3}\left( \frac{n}{6}\right)  \label{w4} \\
&&+\frac{49}{6000}\sigma _{3}\left( \frac{n}{7}\right) +\frac{49}{1500}%
\sigma _{3}\left( \frac{n}{14}\right) +\frac{147}{2000}\sigma _{3}\left(
\frac{n}{21}\right)  \notag \\
&&+\frac{147}{500}\sigma _{3}\left( \frac{n}{42}\right) +(\frac{1}{24}-\frac{%
n}{28})\sigma \left( \frac{n}{6}\right) +(\frac{1}{24}-\frac{n}{24})\sigma
\left( \frac{n}{7}\right)  \notag \\
&&+\frac{1}{70}c_{1}(n)-\frac{407}{2800}c_{2}(n)+\frac{1}{28}c_{3}(n)+\frac{3%
}{140}c_{4}(n)-\frac{1507}{42000}c_{5}(n)  \notag \\
&&-\frac{53}{40}c_{6}(n)+\frac{6899}{1400}c_{7}(n)-\frac{6349}{6000}c_{8}(n)-%
\frac{41}{120}c_{9}(n)  \notag \\
&&+\frac{221}{350}c_{10}(n)-\frac{379}{40}c_{11}(n)+\frac{969}{700}c_{12}(n)-%
\frac{11}{20}c_{13}(n)  \notag \\
&&+\frac{18}{5}c_{14}(n)+\frac{1}{10}c_{15}(n)-\frac{39}{20}c_{16}(n)+\frac{9%
}{28}c_{17}(n)-\frac{3}{10}c_{18}(n)  \notag \\
&&-\frac{4113}{2800}c_{19}(n)+\frac{9}{20}c_{20}(n)  \notag
\end{eqnarray}
\end{theorem}

\begin{proof}
We prove the formula for only $W_{1,42}(n)$ as the rest can be proven
similarly. Glaisher \cite{Glaisher} has proved the following identity%
\begin{equation}
L^{2}(q)=1+\mathop{\displaystyle \sum }\limits_{n=1}^{\infty }\left(
240\sigma _{3}(n)-288n\sigma (n)\right) q^{n}.  \label{lq}
\end{equation}
Replacing $q$ by $q^{42}$ in (\ref{lq}) we have
\begin{equation}
L^{2}(q^{42})=1+\mathop{\displaystyle \sum }\limits_{n=1}^{\infty }\left(
240\sigma _{3}(\frac{n}{42})-\frac{48}{7}n\sigma (\frac{n}{42})\right) q^{n}.
\label{lq2}
\end{equation}
By (\ref{2.1}) we have%
\begin{eqnarray}
L(q)L(q^{42}) &=&\left( 1-24\mathop{\displaystyle \sum }\limits_{n=1}^{%
\infty }\sigma \left( n\right) q^{n}\right) \left( 1-24\mathop{\displaystyle
\sum }\limits_{n=1}^{\infty }\sigma \left( n\right) q^{42n}\right)  \notag \\
&=&1-24\mathop{\displaystyle \sum }\limits_{n=1}^{\infty }\sigma \left(
n\right) q^{n}-24\mathop{\displaystyle \sum }\limits_{n=1}^{\infty }\sigma
\left( \frac{n}{42}\right) q^{n}  \notag \\
&&+576\mathop{\displaystyle \sum }\limits_{n=1}^{\infty }W_{1,42}(n)q^{n}
\label{lq3}
\end{eqnarray}
From (\ref{lq})-(\ref{lq3}) we obtain%
\begin{eqnarray}
(L(q)-42L(q^{42}))^{2} &=&L^{2}(q)-84L(q)L(q^{42})+1764L^{2}(q^{42})
\label{lq4} \\
&=&1681+\mathop{\displaystyle \sum }\limits_{n=1}^{\infty }(240\sigma
_{3}(n)+423360\sigma _{3}(\frac{n}{42})  \notag \\
&&+48384(\frac{1}{24}-\frac{n}{168})\sigma \left( n\right) +48384(\frac{1}{24%
}-\frac{n}{4})\sigma \left( \frac{n}{42}\right)  \notag \\
&&-48384W_{1,42}(n))q^{n}  \notag
\end{eqnarray}
Equating the coefficients of $q^{n}$ on the right hand sides of first part
of Theorem \ref{thm:6} and (\ref{lq4}) we obtain%
\begin{eqnarray*}
&&\frac{28992}{125}\sigma _{3}(n)-\frac{4032}{125}\sigma _{3}(\frac{n}{2})-%
\frac{9072}{125}\sigma _{3}(\frac{n}{3})-\frac{36288}{125}\sigma _{3}(\frac{n%
}{6}) \\
&&-\frac{49392}{125}\sigma _{3}(\frac{n}{7})-\frac{197568}{125}\sigma _{3}(%
\frac{n}{14})-\frac{444528}{125}\sigma _{3}(\frac{n}{21}) \\
&&+\frac{51141888}{125}\sigma _{3}(\frac{n}{42})+\frac{6912}{5}c_{1}(n)+%
\frac{24624}{25}c_{2}(n)+\frac{1728}{5}c_{4}(n) \\
&&+\frac{1008}{125}c_{5}(n)-\frac{296352}{5}c_{6}(n)+\frac{5340384}{25}%
c_{7}(n)-\frac{2067408}{125}c_{8}(n) \\
&&-\frac{34272}{5}c_{9}(n)+\frac{2346624}{25}c_{10}(n)-\frac{3284064}{5}%
c_{11}(n)+\frac{1024128}{25}c_{12}(n) \\
&&-\frac{124992}{5}c_{13}(n)+\frac{653184}{5}c_{14}(n)+\frac{48384}{5}%
c_{15}(n)-\frac{520128}{5}c_{16}(n) \\
&&+\frac{217728}{5}c_{18}(n)+\frac{27216}{25}c_{19}(n)+\frac{36288}{5}%
c_{20}(n) \\
&=&240\sigma _{3}(n)+423360\sigma _{3}(\frac{n}{42})+48384(\frac{1}{24}-%
\frac{n}{168})\sigma \left( n\right) \\
&&+48384(\frac{1}{24}-\frac{n}{4})\sigma \left( \frac{n}{42}\right)
-48384W_{1,42}(n).
\end{eqnarray*}
Solving this equation for $W_{1,42}(n)$ we obtain the asserted formula.
\end{proof}

\section{The Representation Number Formula for $%
x_{1}^{2}+x_{1}x_{2}+x_{2}^{2}+x_{3}^{2}+x_{3}x_{4}+x_{4}^{2}+14(x_{5}^{2}+x_{5}x_{6}+x_{6}^{2}+x_{7}^{2}+x_{7}x_{8}+x_{8}^{2})
$}

\begin{theorem}
\label{thm:8} Let $n\in
\mathbb{N}
$, then%
\begin{eqnarray}
N_{14}(n) &=&\frac{12}{125}\sigma _{3}(n)+\frac{48}{125}\sigma _{3}(\frac{n}{%
2})+\frac{108}{125}\sigma _{3}(\frac{n}{3})+\frac{432}{125}\sigma _{3}(\frac{%
n}{6})  \label{N14} \\
&&+\frac{588}{125}\sigma _{3}(\frac{n}{7})+\frac{2352}{125}\sigma _{3}(\frac{%
n}{14})+\frac{5292}{125}\sigma _{3}(\frac{n}{21})+\frac{21168}{125}\sigma
_{3}(\frac{n}{42})  \notag \\
&&+\frac{288}{25}c_{2}(n)-\frac{6}{25}c_{3}(n)+\frac{294}{25}c_{4}(n)+\frac{%
2592}{175}c_{2}(\frac{n}{3})-\frac{54}{25}c_{3}(\frac{n}{3})  \notag \\
&&-\frac{5778}{175}c_{4}(\frac{n}{3})+\frac{18}{125}c_{5}(n)-\frac{648}{5}%
c_{6}(n)+\frac{2484}{25}c_{7}(n)+\frac{882}{125}c_{8}(n)  \notag \\
&&-\frac{1296}{25}c_{10}(n)-\frac{252}{25}c_{12}(n)+\frac{972}{7}c_{17}(n)+%
\frac{15552}{175}c_{19}(n)  \notag
\end{eqnarray}
\end{theorem}

\begin{proof}
For $l\in \mathbf{%
\mathbb{N}
}_{\mathbf{0}}$ we set,
\begin{equation}
r(l)=\text{card}\left\{ \left( x_{1},...,x_{4}\right) \in
\mathbb{Z}
^{4}:l=x_{1}^{2}+x_{1}x_{2}+x_{2}^{2}+x_{3}^{2}+x_{3}x_{4}+x_{4}^{2}\right\}
,  \label{rl}
\end{equation}
It is known (see for example \cite[Theorem 13]{Huard}) that
\begin{equation}
r(l)=12\sigma (l)-36\sigma (\frac{l}{3}),\text{ }l\in \mathbf{%
\mathbb{N}
}.  \label{r12}
\end{equation}
It is clear from (\ref{1.2}), (\ref{rl}) and (\ref{r12}) that%
\begin{eqnarray*}
N_{14}(n) &=&\mathop{\displaystyle \sum }\limits_{\substack{ l,m\in \mathbf{%
\mathbb{N}
}_{\mathbf{0}}  \\ l+14m=n}}r(l)r(m) \\
&=&r(n)r(0)+r(0)r(\frac{n}{14})+\mathop{\displaystyle \sum }\limits
_{\substack{ l,m\in \mathbf{%
\mathbb{N}
}  \\ l+14m=n}}r(l)(m) \\
&=&12\sigma (n)-36\sigma (\frac{n}{3})+12\sigma (\frac{n}{14})-36\sigma (%
\frac{n}{42}) \\
&&+\mathop{\displaystyle \sum }\limits_{\substack{ l,m\in \mathbf{%
\mathbb{N}
}  \\ l+14m=n}}(12\sigma (l)-36\sigma (\frac{l}{3}))(12\sigma (m)-36\sigma (%
\frac{m}{3}) \\
&=&12\sigma (n)-36\sigma (\frac{n}{3})+12\sigma (\frac{n}{14})-36\sigma (%
\frac{n}{42}) \\
&&+144\mathop{\displaystyle \sum }\limits_{\substack{ l,m\in \mathbf{%
\mathbb{N}
}  \\ l+14m=n}}\sigma (l)\sigma (m)-432\mathop{\displaystyle \sum }\limits
_{\substack{ l,m\in \mathbf{%
\mathbb{N}
}  \\ l+14m=n}}\sigma (\frac{l}{3})\sigma (m) \\
&&-432\mathop{\displaystyle \sum }\limits_{\substack{ l,m\in \mathbf{%
\mathbb{N}
}  \\ l+14m=n}}\sigma (l)\sigma (\frac{m}{3})+1296\mathop{\displaystyle \sum
}\limits_{\substack{ l,m\in \mathbf{%
\mathbb{N}
}  \\ l+14m=n}}\sigma (\frac{l}{3})\sigma (\frac{m}{3}) \\
&=&(12\sigma (n)-36\sigma (\frac{n}{3})+12\sigma (\frac{n}{14})-36\sigma (%
\frac{n}{42})+144W_{1,14}(n) \\
&&-432W_{3,14}(n)-432W_{1,42}(n)+1296W_{1,14}(\frac{n}{3})
\end{eqnarray*}

Formulae for $W_{3,14}(n)$ and $W_{1,42}(n)$ are obtained in this article. $%
W_{1,14}(n)$ was evaluated in \cite{Alaca10}. We also check that it can be
given as%
\begin{eqnarray}
W_{1,14}(n) &=&\frac{1}{600}\sigma _{3}\left( n\right) +\frac{1}{150}\sigma
_{3}\left( \frac{n}{2}\right) +\frac{49}{600}\sigma _{3}\left( \frac{n}{7}%
\right) +\frac{49}{150}\sigma _{3}\left( \frac{n}{14}\right)  \notag \\
&&+(\frac{1}{24}-\frac{n}{56})\sigma \left( n\right) +(\frac{1}{24}-\frac{n}{%
4})\sigma \left( \frac{n}{14}\right) +\frac{2}{175}c_{2}(n)  \notag \\
&&-\frac{1}{600}c_{3}(n)-\frac{107}{4200}c_{4}(n)).  \label{w5}
\end{eqnarray}
Using (\ref{w1}), (\ref{w3}) and (\ref{w5}) we obtained the desired formula.
\end{proof}

For simplicity, taking $u(n)$ as follows
\begin{eqnarray*}
u(n) &=&1680c_{2}(n)+2160c_{2}(\frac{n}{3})-35c_{3}(n)-315c_{3}(\frac{n}{3}%
)+1715c_{4}(n) \\
&&-4815c_{4}(\frac{n}{3})+21c_{5}(n)-18900c_{6}(n)+14490c_{7}(n)+1029c_{8}(n)
\\
&&-7560c_{10}(n)-1470c_{12}+20250c_{17}(n)+12960c_{19}(n)
\end{eqnarray*}%
we may write
\begin{eqnarray}
N_{14}(n) &=&\frac{12}{125}\sigma _{3}(n)+\frac{48}{125}\sigma _{3}(\frac{n}{%
2})+\frac{108}{125}\sigma _{3}(\frac{n}{3})+\frac{432}{125}\sigma _{3}(\frac{%
n}{6})+\frac{588}{125}\sigma _{3}(\frac{n}{7})  \notag \\
&&+\frac{2352}{125}\sigma _{3}(\frac{n}{14})+\frac{5292}{125}\sigma _{3}(%
\frac{n}{21})+\frac{21168}{125}\sigma _{3}(\frac{n}{42})+\frac{6}{875}u(n).
\label{N142}
\end{eqnarray}%
We checked our result for some values of $n$ by using Pari GP.
Denoting the right hand side of (\ref{N142}) by $S_{14}(n)$ we give
the first ten values of $N_{14}(n)$ and $S_{14}(n)$ in Table
\ref{Tab:Tb2} to illustrate the equations.
\begin{table}[th] \caption{A few values of
$N_{14}(n)$.} \label{Tab:Tb2}\centering

\begin{tabular}{|c|c|c|c|c|}
\hline
$n$ & $N_{14}(n)$ & $\sigma _{3}(n)$ & $u(n)$ & $S_{14}(n)$ \\ \hline\hline
$1$ & $12$ & $1$ & $1736$ & $12$ \\
$2$ & $36$ & $9$ & $5068$ & $36$ \\
$3$ & $12$ & $28$ & $1232$ & $12$ \\
$4$ & $84$ & $73$ & $10724$ & $84$ \\
$5$ & $72$ & $126$ & $8736$ & $72$ \\
$6$ & $36$ & $252$ & $-1484$ & $36$ \\
$7$ & $96$ & $344$ & $8498$ & $96$ \\
$8$ & $180$ & $585$ & $13972$ & $180$ \\
$9$ & $12$ & $757$ & $-12376$ & $12$ \\
$10$ & $216$ & $1134$ & $8568$ & $216$ \\ \hline
\end{tabular}%
\end{table}

\section{Acknowledgement}
This material is based upon work supported by the Simons Foundation
Institute Grant Award ID 507536 while the author was in residence at
the Institute for Computational and Experimental Research in
Mathematics in Providence, RI.


\begin{thebibliography}{99}
\bibitem{Glaisher} J. W. L. Glaisher, On the square of the series in which
the coefficients are the sum of the divisors of the exponents, Messenger
Math., 14(1862), 156-163.
\end{thebibliography}
\end{document}